# IMAGINARY QUADRATIC FIELDS $k$ WITH $\mathrm{Cl}_2(k) \simeq (2, 2^m)$ AND RANK $\mathrm{Cl}_2(k^1) = 2$

E. BENJAMIN, F. LEMMERMEYER, C. SNYDER

ABSTRACT. Let $k$ be an imaginary quadratic number field and $k^1$ the Hilbert 2-class field of $k$. We give a characterization of those $k$ with $\mathrm{Cl}_2(k) \simeq (2, 2^m)$ such that $\mathrm{Cl}_2(k^1)$ has 2 generators.

## 1. INTRODUCTION

Let $k$ be an algebraic number field with $\mathrm{Cl}_2(k)$, the Sylow 2-subgroup of its ideal class group, $\mathrm{Cl}(k)$. Denote by $k^1$ the Hilbert 2-class field of $k$ (in the wide sense). Also let $k^n$ (for $n$ a nonnegative integer) be defined inductively as: $k^0 = k$ and $k^{n+1} = (k^n)^1$. Then

$$k^0 \subseteq k^1 \subseteq k^2 \subseteq \cdots \subseteq k^n \subseteq \cdots$$

is called the 2-class field tower of $k$. If $n$ is the minimal integer such that $k^n = k^{n+1}$, then $n$ is called the length of the tower. If no such $n$ exists, then the tower is said to be of infinite length.

At present there is no known decision procedure to determine whether or not the (2-)class field tower of a given field $k$ is infinite. However, it is known by group theoretic results (see [2]) that if $\operatorname{rank} \mathrm{Cl}_2(k^1) \leq 2$, then the tower is finite, in fact of length at most 3. (Here the rank means minimal number of generators.) On the other hand, until now (see Table 1 and the penultimate paragraph of this introduction) all examples in the mathematical literature of imaginary quadratic fields with $\operatorname{rank} \mathrm{Cl}_2(k^1) \geq 3$ (let us mention in particular Schmithals [13]) have infinite 2-class field tower. Nevertheless, if we are interested in developing a decision procedure for determining if the 2-class field tower of a field is infinite, then a good starting point would be to find a procedure for sieving out those fields with $\operatorname{rank} \mathrm{Cl}_2(k^1) \leq 2$. We have already started this program for imaginary quadratic number fields $k$. In [1] we classified all imaginary quadratic fields whose 2-class field $k^1$ has cyclic 2-class group. In this paper we determine when $\mathrm{Cl}_2(k^1)$ has rank 2 for imaginary quadratic fields $k$ with $\mathrm{Cl}_2(k)$ of type $(2, 2^m)$. (The notation $(2, 2^m)$ means the direct sum of a group of order 2 and a cyclic group of order $2^m$.) The group theoretic results mentioned above also show that such fields have 2-class field tower of length 2.

From a classification of imaginary quadratic number fields $k$ with $\mathrm{Cl}_2(k) \simeq (2, 2^m)$ and our results from [1] we see that it suffices to consider discriminants $d = d_1 d_2 d_3$ with prime discriminants $d_1, d_2 > 0$, $d_3 < 0$ such that exactly one of the $(d_i/p_j)$ equals $-1$ (we let $p_j$ denote the prime dividing $d_j$); thus there are only two cases:

A) $(d_1/p_2) = (d_1/p_3) = +1$, $(d_2/p_3) = -1$;

  B) $(d_1/p_3) = (d_2/p_3) = +1$, $(d_1/p_2) = -1$.

The $C_4$-factorization corresponding to the nontrivial 4-part of $\mathrm{Cl}_2(k)$ is $d = d_1 \cdot d_2 d_3$ in case A) and $d = d_1 d_2 \cdot d_3$ in case B). Note that, by our results from [1], some of these fields have cyclic $\mathrm{Cl}_2(k^1)$; however, we do not exclude them right from the start since there is no extra work involved and since it provides a welcome check on our earlier work.

The main result of the paper is that $\operatorname{rank} \mathrm{Cl}_2(k^1) = 2$ only occurs for fields of type B); more precisely, we prove the following

**Theorem 1.** *Let $k$ be a complex quadratic number field with $\mathrm{Cl}_2(k) \simeq (2, 2^m)$, and let $k^1$ be its $2$-class field. Then $\operatorname{rank} \mathrm{Cl}_2(k^1) = 2$ if and only if $\operatorname{disc} k = d_1 d_2 d_3$ is the product of three prime discriminants $d_1, d_2 > 0$ and $-4 \neq d_3 < 0$ such that $(d_1/p_3) = (d_2/p_3) = +1$, $(d_1/p_2) = -1$, and $h_2(K) = 2$, where $K$ is a nonnormal quartic subfield of one of the two unramified cyclic quartic extensions of $k$ such that $\mathbb{Q}(\sqrt{d_1 d_2}) \subset K$.*

This result is the first step in the classification of imaginary quadratic number fields $k$ with $\operatorname{rank} \mathrm{Cl}_2(k^1) = 2$; it remains to solve these problems for fields with $\operatorname{rank} \mathrm{Cl}_2(k) = 3$ and those with $\mathrm{Cl}_2(k) \supseteq (4, 4)$ since we know that $\operatorname{rank} \mathrm{Cl}_2(k^1) \geq 5$ whenever $\operatorname{rank} \mathrm{Cl}_2(k) \geq 4$ (using Schur multipliers as in [1]).

As a demonstration of the utility of our results, we give in Table 1 below a list of the first 12 imaginary quadratic fields $k$, arranged by decreasing value of their discriminants, with $\operatorname{rank} \mathrm{Cl}_2(k) = 2$ and noncyclic $\mathrm{Cl}_2(k^1)$.

TABLE 1

| disc $k$ | factors | $\mathrm{Cl}_2(k)$ | type | $f$ | $\mathrm{Cl}_2(K)$ | $r$ | $\mathrm{Cl}_2(k_{\mathrm{gen}})$ |
|---|---|---|---|---|---|---|---|
| $-1015$ | $-7 \cdot 5 \cdot 29$ | $(2,8)$ | A | $x^4 - 22x^2 + 261$ | $(4)$ | $\geq 3$ | $(2,2,8)$ |
| $-1240$ | $-31 \cdot 8 \cdot 5$ | $(2,4)$ | B | $x^4 - 6x^2 - 31$ | $(2)$ | $2$ | $(2,2,8)$ |
| $-1443$ | $-3 \cdot 13 \cdot 37$ | $(2,4)$ | B | $x^4 - 86x^2 - 75$ | $(2)$ | $2$ | $(2,2,8)$ |
| $-1595$ | $-11 \cdot 5 \cdot 29$ | $(2,8)$ | A | $x^4 + 26x^2 + 1445$ | $(4)$ | $\geq 3$ | $(2,2,8)$ |
| $-1615$ | $-19 \cdot 5 \cdot 17$ | $(2,4)$ | B | $x^4 + 26x^2 - 171$ | $(2)$ | $2$ | $(2,2,8)$ |
| $-1624$ | $-7 \cdot 8 \cdot 29$ | $(2,8)$ | B | $x^4 - 30x^2 - 7$ | $(2)$ | $2$ | $(2,2,8)$ |
| $-1780$ | $-4 \cdot 5 \cdot 89$ | $(2,4)$ | A | $x^4 + 6x^2 + 89$ | $(4)$ | $3$ | $(2,2,4)$ |
| $-2035$ | $-11 \cdot 5 \cdot 37$ | $(2,4)$ | B | $x^4 - 54x^2 - 11$ | $(2,2)$ | $3$ | $(2,2,16)$ |
| $-2067$ | $-3 \cdot 13 \cdot 53$ | $(2,4)$ | A | $x^4 + x^2 + 637$ | $(2,2)$ | $3$ | $(2,2,4)$ |
| $-2072$ | $-7 \cdot 8 \cdot 37$ | $(2,8)$ | B | $x^4 + 34x^2 - 7$ | $(2,2)$ | $\geq 3$ | $(2,2,8)$ |
| $-2379$ | $-3 \cdot 13 \cdot 61$ | $(4,4)$ | | | | $\geq 3$ | $(4,4,8)$ |
| $-2392$ | $-23 \cdot 8 \cdot 13$ | $(2,4)$ | B | $x^4 + 18x^2 - 23$ | $(2)$ | $2$ | $(2,2,8)$ |

Here $f$ denotes a generating polynomial for a field $K$ as in Theorem 1, $r$ denotes the rank of $\mathrm{Cl}_2(k^1)$. The cases where $r = 3$ follow from our theorem combined with Blackburn's upper bound for the number of generators of derived groups (it implies that finite 2-groups $G$ with $G/G' \simeq (2,4)$ satisfy $\operatorname{rank} G' \leq 3$), see [3].

In order to verify that $\mathrm{Cl}_2(k^1)$ has rank at least 3 for $k = \mathbb{Q}(\sqrt{-2379})$ it is sufficient to show that its genus class field $k_{\mathrm{gen}}$ has class group $(4, 4, 8)$: in fact, $\mathrm{Cl}_2(k^1)$ then contains a quotient of $(4, 4, 8)$ by $(2, 2) \simeq \operatorname{Gal}(k^1/k_{\mathrm{gen}})$, and the claim follows.

We mention one last feature gleaned from the table. It follows from conditional Odlyzko bounds (assuming the Generalized Riemann Hypothesis) that those quadratic fields with rank $\text{Cl}_2(k^1) \geq 3$ and discriminant $0 > d > -2000$ have finite class field tower; unconditional proofs are not known. Hence, conditionally, we conclude that those $k$ with discriminants $-1015, -1595$ and $-1780$ have finite (2-)class field tower even though rank $\text{Cl}_2(k^1) \geq 3$. Of course, it would be interesting to determine the length of their towers.

The structure of this paper is as follows: we use results from group theory developed in Section 2 to pull down the condition rank $\text{Cl}_2(k^1) = 2$ from the field $k^1$ with degree $2^{m+2}$ to a subfield $L$ of $k^1$ with degree 8. Using the arithmetic of dihedral fields from Section 4 we then go down to the field $K$ of degree 4 occurring in Theorem 1.

## 2. Group Theoretic Preliminaries

Let $G$ be a group. If $x, y \in G$, then we let $[x, y] = x^{-1}y^{-1}xy$ denote the commutator of $x$ and $y$. If $A$ and $B$ are nonempty subsets of $G$, then $[A, B]$ denotes the subgroup of $G$ generated by the set $\{[a, b] : a \in A, b \in B\}$. The lower central series $\{G_j\}$ of $G$ is defined inductively by: $G_1 = G$ and $G_{j+1} = [G, G_j]$ for $j \geq 1$. The derived series $\{G^{(n)}\}$ is defined inductively by: $G^{(0)} = G$ and $G^{(n+1)} = [G^{(n)}, G^{(n)}]$ for $n \geq 0$. Notice that $G^{(1)} = G_2 = [G, G]$ the commutator subgroup, $G'$, of $G$.

Throughout this section, we assume that $G$ is a finite, nonmetacyclic, 2-group such that its abelianization $G^{\text{ab}} = G/G'$ is of type $(2, 2^m)$ for some positive integer $m$ (necessarily $\geq 2$). Let $G = \langle a, b \rangle$, where $a^2 \equiv b^{2^m} \equiv 1 \mod G_2$ (actually $\mod G_3$ since $G$ is nonmetacyclic, cf. [1]); $c_2 = [a, b]$ and $c_{j+1} = [b, c_j]$ for $j \geq 2$.

**Lemma 1.** *Let $G$ be as above (but not necessarily metabelian). Suppose that $d(G') = n$ where $d(G')$ denotes the minimal number of generators of the derived group $G' = G_2$ of $G$. Then*
$$G' = \langle c_2, c_3, \cdots, c_{n+1} \rangle;$$
*moreover,*
$$G_2/G_2^2 \simeq \langle c_2 G_2^2 \rangle \oplus \cdots \oplus \langle c_{n+1} G_2^2 \rangle.$$

*Proof.* By the Burnside Basis Theorem, $d(G_2) = d(G_2/\Phi(G))$, where $\Phi(G)$ is the Frattini subgroup of $G$, i.e. the intersection of all maximal subgroups of $G$, see [5]. But in the case of a 2-group, $\Phi(G) = G^2$, see [8]. By Blackburn, [3], since $G/G_2^2$ has elementary derived group, we know that $G_2/G_2^2 \simeq \langle c_2 G_2^2 \rangle \oplus \cdots \oplus \langle c_{n+1} G_2^2 \rangle$. Again, by the Burnside Basis Theorem, $G_2 = \langle c_2, \cdots, c_{n+1} \rangle$. □

**Lemma 2.** *Let $G$ be as above. Moreover, assume $G$ is metabelian. Let $H$ be a maximal subgroup of $G$ such that $H/G'$ is cyclic, and denote the index $(G' : H')$ by $2^\kappa$. Then $G'$ contains an element of order $2^\kappa$.*

*Proof.* Without loss of generality, let $H = \langle b, G' \rangle$. Notice that $G' = \langle c_2, c_3, \cdots \rangle$ and by our presentation of $H$, $H' = \langle c_3, c_4, \cdots \rangle$. Thus, $G'/H' = \langle c_2 H' \rangle$. But since $(G' : H') = 2^\kappa$, the order of $c_2$ is $\geq 2^\kappa$. This establishes the lemma. □

**Lemma 3.** *Let $G$ be as above and again assume $G$ is metabelian. Let $H$ be a maximal subgroup of $G$ such that $H/G'$ is cyclic, and assume that $(G' : H') \equiv 0 \mod 4$. If $d(G') = 2$, then $G_2 = \langle c_2, c_3 \rangle$ and $G_j = \langle c_2^{2^{j-2}}, c_3^{2^{j-3}} \rangle$ for $j > 2$.*

*Proof.* Assume that $d(G') = 2$. By Lemma 1, $G_2 = \langle c_2, c_3 \rangle$ and hence $c_4 \in \langle c_2, c_3 \rangle$. Write $c_4 = c_2^x c_3^y$ where $x, y$ are positive integers. Without loss of generality, let $H = \langle b, c_2, c_3 \rangle$ and write $(G' : H') = 2^\kappa$ for some $\kappa \geq 2$. Since $c_3, c_4 \in H'$ we have, $c_2^x \equiv 1 \mod H'$. By the proof of Lemma 2, this implies that $x \equiv 0 \mod 2^\kappa$. Write $x = 2^\kappa x_1$ for some positive integer $x_1$. On the other hand, since $c_4, c_2^{2^\kappa x_1} \in G_4$, we see that $c_3^y \equiv 1 \mod G_4$. If $y$ were odd, then $c_3 \in G_4$. This, however, implies that $G_2 = \langle c_2 \rangle$, contrary to our assumptions. Thus $y$ is even, say $y = 2y_1$. From all of this we see that $c_4 = c_2^{2^\kappa x_1} c_3^{2y_1}$. Consequently, by induction we have $c_j \in \langle c_2^{2^{j-2}}, c_3^{2^{j-3}} \rangle$ for all $j \geq 4$. Since $G_j = \langle c_2^{2^{j-2}}, c_3^{2^{j-3}}, \cdots, c_{j-1}^2, c_j, c_{j+1}, \cdots \rangle$, cf. [1], we obtain the lemma. □

Let us translate the above into the field-theoretic language. Let $k$ be an imaginary quadratic number field of type A) or B) (see the Introduction), and let $M/k$ be one of the two quadratic subextensions of $k^1/k$ over which $k^1$ is cyclic. If $h_2(M) = 2^{m+\kappa}$ and $\text{Cl}_2(k) = (2, 2^m)$, then Lemma 2 implies that $\text{Cl}_2(k^1)$ contains an element of order $2^\kappa$. Table 2 contains the relevant information for the fields occurring in Table 1. An application of the class number formula to $M/\mathbb{Q}$ (see e.g. Proposition 3 below) shows immediately that $h_2(M) = 2^{m+\kappa}$, where $2^\kappa$ is the class number of the quadratic subfield $\mathbb{Q}(\sqrt{d_i d_j})$ of $M$, where $(d_i/p_j) = +1$; in particular, we always have $\kappa \geq 2$, and the assumption $(G' : H') \geq 4$ is always satisfied for the fields that we consider.

TABLE 2

| $M_1$ | $\text{Cl}_2(M_1)$ | $M_2$ | $\text{Cl}_2(M_2)$ |
|---|---|---|---|
| $\mathbb{Q}(\sqrt{5}, \sqrt{-7 \cdot 29})$ | $(2, 16)$ | $\mathbb{Q}(\sqrt{5 \cdot 29}, \sqrt{-7})$ | $(2, 16)$ |
| $\mathbb{Q}(\sqrt{2}, \sqrt{-5 \cdot 31})$ | $(4, 4)$ | $\mathbb{Q}(\sqrt{5}, \sqrt{-2 \cdot 31})$ | $(2, 16)$ |
| $\mathbb{Q}(\sqrt{13}, \sqrt{-3 \cdot 37})$ | $(2, 16)$ | $\mathbb{Q}(\sqrt{37}, \sqrt{-3 \cdot 13})$ | $(2, 16)$ |
| $\mathbb{Q}(\sqrt{-11}, \sqrt{5 \cdot 29})$ | $(2, 16)$ | $\mathbb{Q}(\sqrt{29}, \sqrt{-5 \cdot 11})$ | $(2, 16)$ |
| $\mathbb{Q}(\sqrt{5}, \sqrt{-17 \cdot 19})$ | $(4, 4)$ | $\mathbb{Q}(\sqrt{17}, \sqrt{-5 \cdot 19})$ | $(2, 16)$ |
| $\mathbb{Q}(\sqrt{29}, \sqrt{-2 \cdot 7})$ | $(2, 16)$ | $\mathbb{Q}(\sqrt{2}, \sqrt{-7 \cdot 29})$ | $(2, 16)$ |
| $\mathbb{Q}(\sqrt{5 \cdot 89}, \sqrt{-1})$ | $(4, 4)$ | $\mathbb{Q}(\sqrt{5}, \sqrt{-89})$ | $(2, 8)$ |
| $\mathbb{Q}(\sqrt{37}, \sqrt{-5 \cdot 11})$ | $(4, 4)$ | $\mathbb{Q}(\sqrt{5}, \sqrt{-37 \cdot 11})$ | $(2, 32)$ |
| $\mathbb{Q}(\sqrt{53}, \sqrt{-3 \cdot 13})$ | $(4, 4)$ | $\mathbb{Q}(\sqrt{13 \cdot 53}, \sqrt{-3})$ | $(2, 2, 4)$ |
| $\mathbb{Q}(\sqrt{37}, \sqrt{-2 \cdot 7})$ | $(2, 16)$ | $\mathbb{Q}(\sqrt{2}, \sqrt{-7 \cdot 37})$ | $(2, 16)$ |
| $\mathbb{Q}(\sqrt{13}, \sqrt{-2 \cdot 23})$ | $(4, 4)$ | $\mathbb{Q}(\sqrt{2}, \sqrt{-13 \cdot 23})$ | $(2, 16)$ |

We now use the above results to prove the following useful proposition.

**Proposition 1.** *Let $G$ be a nonmetacyclic 2-group such that $G/G' \simeq (2, 2^m)$; (hence $m > 1$). Let $H$ and $K$ be the two maximal subgroups of $G$ such that $H/G'$ and $K/G'$ are cyclic. Moreover, assume that $(G' : H') \equiv 0 \mod 4$. Finally, assume* *that $N$ is a subgroup of index $4$ in $G$ not contained in $H$ or $K$ Then*

$$(N : N') \begin{cases} = & 2^m & \text{if } d(G') = 1 \\ = & 2^{m+1} & \text{if } d(G') = 2 \\ \geq & 2^{m+2} & \text{if } d(G') \geq 3 \end{cases}.$$

*Proof.* Without loss of generality we assume that $G$ is metabelian. Let $G = \langle a, b \rangle$, where $a^2 \equiv b^{2^m} \equiv 1 \mod G_3$. Also let $H = \langle b, G' \rangle$ and $K = \langle ab, G' \rangle$ (without loss of generality). Then $N = \langle ab^2, G' \rangle$ or $N = \langle a, b^4, G' \rangle$.

Suppose that $N = \langle ab^2, G' \rangle$.

First assume $d(G') = 1$. Then $G' = \langle c_2 \rangle$ and thus $N' = \langle [ab^2, c_2] \rangle$. But $[ab^2, c_2] = c_2^2 \eta_4$ for some $\eta_4 \in G_4 = \langle c_2^4 \rangle$ (cf. Lemma 1 of [1]). Hence, $N' = \langle c_2^2 \rangle$, and so $(G' : N') = 2$. Since $(N : G') = 2^{m-1}$, we get $(N : N') = 2^m$ as desired.

Next, assume that $d(G') = 2$. Then $N = \langle ab^2, c_2, c_3 \rangle$ by Lemma 1. Notice that $[ab^2, c_2] = c_2^2 \eta_4$ and $[ab^2, c_3] = c_3^2 \eta_5$ where $\eta_j \in G_j$ for $j = 4, 5$. Hence $N' = \langle c_2^2 \eta_4, c_3^2 \eta_5, N_3 \rangle$ and so $\langle c_2^2 \eta_4, c_3^2 \eta_5 \rangle \subseteq N'$. But then $N'G_5 \supseteq \langle c_2^4, c_3^2 \rangle = G_4$ by Lemma 3. Therefore, by [5], $N' \supseteq G_4$. But notice that $N_3 \subseteq G_4$. Thus $N' = \langle c_2^2, c_3^2 \rangle$ and so $(G' : N') = 4$ which in turn implies that $(N : N') = 2^{m+1}$, as desired.

Finally, assume $d(G') \geq 3$. Then $d(G'/G_5) = 3$. Moreover there exists an exact sequence

$$N/N' \longrightarrow (N/G_5)/(N/G_5)' \longrightarrow 1,$$

and thus $\#N^{\text{ab}} \geq \#(N/G_5)^{\text{ab}}$. Hence it suffices to prove the result for $G_5 = 1$ which we now assume. $N = \langle ab^2, c_2, c_3, c_4 \rangle$ and so, arguing as above, we have $N' = \langle c_2^2 \eta_4, c_3^2 \eta_5, c_4^2 \eta_6, N_3 \rangle = \langle c_2^2 \eta_4, c_3^2, N_3 \rangle$, where $\eta_j \in G_j$. But $N_3 = \langle [ab^2, c_2^2 \eta_4] \rangle = \langle c_2^4 \rangle$. Therefore, $N' = \langle c_2^2 \eta_4, c_3^2 \rangle$. From this we see that $(G' : N') = 8$ and thus $(N : N') = 2^{m+2}$ as desired.

Now suppose that $N = \langle a, b^4, G' \rangle$. Then the proof is essentially the same as above once we notice that $[a, b^4] \equiv c_3^2 c_2^{-4} \mod G_5$.

This establishes the proposition. $\square$

## 3. Number Theoretic Preliminaries

**Proposition 2.** *Let $K/k$ be a quadratic extension, and assume that the class number of $k$, $h(k)$, is odd. If $K$ has an unramified cyclic extension $M$ of order $4$, then $M/k$ is normal and $\operatorname{Gal}(M/k) \simeq D_4$.*

*Proof.* Rédei and Reichardt [12] proved this for $k = \mathbb{Q}$; the general case is analogous. $\square$

We shall make extensive use of the class number formula for extensions of type $(2, 2)$:

**Proposition 3.** *Let $K/k$ be a normal quartic extension with Galois group of type $(2,2)$, and let $k_j$ $(j = 1, 2, 3)$ denote the quadratic subextensions. Then*

(1) $$h(K) = 2^{d-\kappa-2-\upsilon} q(K) h(k_1) h(k_2) h(k_3)/h(k)^2,$$

*where $q(K) = (E_K : E_1 E_2 E_3)$ denotes the unit index of $K/k$ ($E_j$ is the unit group of $k_j$), $d$ is the number of infinite primes in $k$ that ramify in $K/k$, $\kappa$ is the $\mathbb{Z}$-rank of the unit group $E_k$ of $k$, and $\upsilon = 0$ except when $K \subseteq k(\sqrt{E_k})$, where $\upsilon = 1$.*

*Proof.* See [10]. $\square$

Another important result is the ambiguous class number formula. For cyclic extensions $K/k$, let $\mathrm{Am}(K/k)$ denote the group of ideal classes in $K$ fixed by $\mathrm{Gal}(K/k)$, i.e. the ambiguous ideal class group of $K$, and $\mathrm{Am}_2$ its 2-Sylow subgroup.

**Proposition 4.** *Let $K/k$ be a cyclic extension of prime degree $p$; then the number of ambiguous ideal classes is given by*

$$\# \mathrm{Am}(K/k) = h(k) \frac{p^{t-1}}{(E:H)},$$

*where $t$ is the number of primes (including those at $\infty$) of $k$ that ramify in $K/k$, $E$ is the unit group of $k$, and $H$ is its subgroup consisting of norms of elements from $K^\times$. Moreover, $\mathrm{Cl}_p(K)$ is trivial if and only if $p \nmid \# \mathrm{Am}(K/k)$.*

*Proof.* See Lang [9, part II] for the formula. For a proof of the second assertion (see e.g. Moriya [11]), note that $\mathrm{Am}(K/k)$ is defined by the exact sequence

$$1 \longrightarrow \mathrm{Am}(K/k) \longrightarrow \mathrm{Cl}(K) \longrightarrow \mathrm{Cl}(K)^{1-\sigma} \longrightarrow 1,$$

where $\sigma$ generates $\mathrm{Gal}(K/k)$. Taking $p$-parts we see that $p \nmid \# \mathrm{Am}(K/k)$ is equivalent to $\mathrm{Cl}_p(K) = \mathrm{Cl}_p(K)^{1-\sigma}$. By induction we get $\mathrm{Cl}_p(K) = \mathrm{Cl}_p(K)^{(1-\sigma)^p}$, but since $(1-\sigma)^p \equiv 0 \bmod p$ in the group ring $\mathbb{Z}[G]$, this implies $\mathrm{Cl}_p(K) \subseteq \mathrm{Cl}_p(K)^p$. But then $\mathrm{Cl}_p(K)$ must be trivial. □

We make one further remark concerning the ambiguous class number formula that will be useful below. If the class number $h(k)$ is odd, then it is known that $\# \mathrm{Am}_2(K/k) = 2^r$ where $r = \mathrm{rank}\, \mathrm{Cl}_2(K)$.

We also need a result essentially due to G. Gras [4]:

**Proposition 5.** *Let $K/k$ be a quadratic extension of number fields and assume that $h_2(k) = \# \mathrm{Am}_2(K/k) = 2$. Then $K/k$ is ramified and*

$$\mathrm{Cl}_2(K) \simeq \begin{cases} (2,2) \text{ or } \mathbb{Z}/2^n\mathbb{Z} \ (n \geq 3) & \text{if } \#\kappa_{K/k} = 1, \\ \mathbb{Z}/2^n\mathbb{Z} \ (n \geq 1) & \text{if } \#\kappa_{K/k} = 2, \end{cases}$$

*where $\kappa_{K/k}$ denotes the set of ideal classes of $k$ that become principal (capitulate) in $K$.*

*Proof.* We first notice that $K/k$ is ramified. If the extension were unramified, then $K$ would be the 2-class field of $k$, and since $\mathrm{Cl}_2(k)$ is cyclic, it would follow that $\mathrm{Cl}_2(K) = 1$, contrary to assumption.

Before we start with the rest of the proof, we cite the results of Gras that we need (we could also give a slightly longer direct proof without referring to his results). Let $K/k$ be a cyclic extension of prime power order $p^r$, and let $\sigma$ be a generator of $G = \mathrm{Gal}(K/k)$. For any $p$-group $M$ on which $G$ acts we put $M_i = \{m \in M : m^{(1-\sigma)^i} = 1\}$. Moreover, let $\nu$ be the algebraic norm, that is, exponentiation by $1 + \sigma + \sigma^2 + \ldots + \sigma^{p^r-1}$. Then [4, Cor. 4.3] reads

**Lemma 4.** *Suppose that $M^\nu = 1$; let $n$ be the smallest positive integer such that $M_n = M$ and write $n = a(p-1) + b$ with integers $a \geq 0$ and $0 \leq b \leq p - 2$. If $\# M_{i+1}/M_i = p$ for $i = 0, 1, \ldots, n-1$, then $M \simeq (\mathbb{Z}/p^{a+1}\mathbb{Z})^b \times (\mathbb{Z}/p^a\mathbb{Z})^{p-1-b}$.*

We claim that if $\kappa_{K/k} = 2$, then $M = \mathrm{Cl}_2(K)$ satisfies the assumptions of Lemma 4: in fact, let $j = j_{k \to K}$ denote the transfer of ideal classes. Then $c^{1+\sigma} = j(N_{K/k} c)$ for any ideal class $c \in \mathrm{Cl}_2(K)$, hence $M^\nu = j(\mathrm{Cl}_2(k)) = 1$. Moreover, $M_1 = \operatorname{Am}_2(K/k)$ in our case, hence $M_1/M_0$ has order 2. Since the orders of $M_{i+1}/M_i$ decrease towards 1 as $i$ grows (Gras [4, Prop. 4.1.ii)]), we conclude that $\# M_{i+1}/M_i = 2$ for all $i < n$. Since $a = n$ and $b = 0$ when $p = 2$, Lemma 4 now implies that $\operatorname{Cl}_2(K) \simeq \mathbb{Z}/2^n\mathbb{Z}$, that is, the 2-class group is cyclic.

The second result of Gras that we need is [4, Prop. 4.3]

**Lemma 5.** *Suppose that $M^\nu \neq 1$ but assume the other conditions in Lemma 4. Then $n \geq 2$ and*

$$M \simeq \begin{cases} (\mathbb{Z}/p^2\mathbb{Z}) \times (\mathbb{Z}/p\mathbb{Z})^{n-2} & \text{if } n < p; \\ (\mathbb{Z}/p\mathbb{Z})^p \text{ or } (\mathbb{Z}/p^2\mathbb{Z}) \times (\mathbb{Z}/p\mathbb{Z})^{n-2} & \text{if } n = p; \\ (\mathbb{Z}/p^{a+1}\mathbb{Z})^b \times (\mathbb{Z}/p^a\mathbb{Z})^{p-1-b} & \text{if } n > p. \end{cases}$$

If $\kappa_{K/k} = 1$, then this lemma shows that $\operatorname{Cl}_2(K)$ is either cyclic of order $\geq 4$ or of type $(2,2)$. (Notice that the hypothesis of the lemma is satisfied since $K/k$ is ramified implying that the norm $N_{K/k} : \operatorname{Cl}_2(K) \longrightarrow \operatorname{Cl}_2(k)$ is onto; and so the argument above this lemma applies.) It remains to show that the case $\operatorname{Cl}_2(K) \simeq \mathbb{Z}/4\mathbb{Z}$ cannot occur here.

Now assume that $\operatorname{Cl}_2(K) = \langle C \rangle \simeq \mathbb{Z}/4\mathbb{Z}$; since $K/k$ is ramified, the norm $N_{K/k} : \operatorname{Cl}_2(K) \longrightarrow \operatorname{Cl}_2(k)$ is onto, and using $\kappa_{K/k} = 1$ once more we find $C^{1+\sigma} = c$, where $c$ is the nontrivial ideal class from $\operatorname{Cl}_2(k)$. On the other hand, $c \in \operatorname{Cl}_2(k)$ still has order 2 in $\operatorname{Cl}_2(K)$, hence we must also have $C^2 = C^{1+\sigma}$. But this implies that $C^\sigma = C$, i.e. that each ideal class in $K$ is ambiguous, contradicting our assumption that $\# \operatorname{Am}_2(K/k) = 2$. $\square$

## 4. Arithmetic of some Dihedral Extensions

In this section we study the arithmetic of some dihedral extensions $L/\mathbb{Q}$, that is, normal extensions $L$ of $\mathbb{Q}$ with Galois group $\operatorname{Gal}(L/\mathbb{Q}) \simeq D_4$, the dihedral group of order 8. Hence $D_4$ may be presented as $\langle \tau, \sigma | \tau^2 = \sigma^4 = 1, \tau\sigma\tau = \sigma^{-1} \rangle$. Now consider the following diagrams (Galois correspondence):

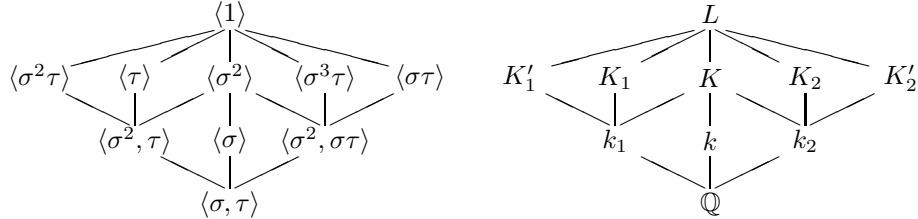

In this situation, we let $q_1 = (E_L : E_1 E_1' E_K)$ and $q_2 = (E_L : E_2 E_2' E_K)$ denote the unit indices of the bicyclic extensions $L/k_1$ and $L/k_2$, where $E_i$ and $E_i'$ are the unit groups in $K_i$ and $K_i'$ respectively. Finally, let $\kappa_i$ denote the kernel of the transfer of ideal classes $j_{k_i \to K_i} : \operatorname{Cl}_2(k_i) \longrightarrow \operatorname{Cl}_2(K_i)$ for $i = 1, 2$.

The following remark will be used several times: if $K_1 = k_1(\sqrt{\alpha})$ for some $\alpha \in k_1$, then $k_2 = \mathbb{Q}(\sqrt{a})$, where $a = \alpha\alpha'$ is the norm of $\alpha$. To see this, let $\gamma = \sqrt{\alpha}$; then $\gamma^\tau = \gamma$, since $\gamma \in K_1$. Clearly $\gamma^{1+\sigma} = \sqrt{a} \in K$ and hence fixed by $\sigma^2$. Furthermore,

$$(\gamma^{1+\sigma})^{\sigma\tau} = \gamma^{\sigma\tau+\sigma^2\tau} = \gamma^{\tau\sigma^3+\tau\sigma^2} = (\gamma^\tau)^{\sigma^3+\sigma^2} = \gamma^{\sigma^3+\sigma^2} = \gamma^{(1+\sigma)\sigma^2} = \gamma^{1+\sigma},$$

implying that $\sqrt{a} \in k_2$. Finally notice that $\sqrt{a} \notin \mathbb{Q}$, since otherwise $\sqrt{\alpha'} = \sqrt{a}/\sqrt{\alpha} \in K_1$ implying that $K_1/\mathbb{Q}$ is normal, which is not the case.

Recall that a quadratic extension $K = k(\sqrt{\alpha}\,)$ is called *essentially ramified* if $\alpha \mathcal{O}_k$ is not an ideal square. This definition is independent of the choice of $\alpha$.

**Proposition 6.** *Let $L/\mathbb{Q}$ be a non-CM totally complex dihedral extension not containing $\sqrt{-1}$, and assume that $L/K_1$ and $L/K_2$ are essentially ramified. If the fundamental unit of the real quadratic subfield of $K$ has norm $-1$, then $q_1 q_2 = 2$.*

*Proof.* Notice first that $k$ cannot be real (in fact, $K$ is not totally real by assumption, and since $L/k$ is a cyclic quartic extension, no infinite prime can ramify in $K/k$); thus exactly one of $k_1$, $k_2$ is real, and the other is complex. Multiplying the class number formulas, Proposition 3, for $L/k_1$ and $L/k_2$ (note that $\upsilon = 0$ since both $L/K_1$ and $L/K_2$ are essentially ramified) we find that $2q_1 q_2$ is a square. If we can prove that $q_1, q_2 \leq 2$, then $2q_1 q_2$ is a square between 2 and 8, which implies that we must have $2q_1 q_2 = 4$ and $q_1 q_2 = 2$ as claimed.

We start by remarking that if $\zeta \eta$ becomes a square in $L$, where $\zeta$ is a root of unity in $L$, then so does one of $\pm \eta$. This follows from the fact that the only non-trivial roots of unity that can be in $L$ are the sixth roots of unity $\langle \zeta_6 \rangle$, and here $\zeta_6 = -\zeta_3^2$.

Now we prove that $q_1 \leq 2$ under the assumptions we made; the claim $q_2 \leq 2$ will then follow by symmetry. Assume first that $k_1$ is real and let $\varepsilon$ be the fundamental unit of $k_1$. We claim that $\sqrt{\pm \varepsilon} \notin L$. Suppose otherwise; then $k_1(\sqrt{\pm \varepsilon}\,)$ is one of $K_1$, $K_1'$ or $K$. If $k_1(\sqrt{\pm \varepsilon}\,) = K_1$, then $K_1' = k_1(\sqrt{\pm \varepsilon'}\,)$ and $K = k_1(\sqrt{\varepsilon \varepsilon'}\,)$. (Here and below $x' = x^\sigma$.) This however cannot occur since by assumption $\varepsilon \varepsilon' = -1$ implying that $\sqrt{-1} \in L$, a contradiction. Similarly, if $k_1(\sqrt{\pm \varepsilon}\,) = K$, then again $\sqrt{-1} \in L$.

Thus $\sqrt{\pm \varepsilon} \notin L$, and $E_1 = \langle -1, \varepsilon, \eta \rangle$ for some unit $\eta \in E_1$. Suppose that $\sqrt{u\eta} \in L$ for some unit $u \in k_1$. Then $L = K_1(\sqrt{u\eta}\,)$, contradicting our assumption that $L/K_1$ is essentially ramified. The same argument shows that $\sqrt{u\eta'} \notin L$, hence either $E_L = \langle \zeta, \varepsilon, \eta, \eta' \rangle$ and $q_1 = 1$ or $E_L = \langle \zeta, \varepsilon, \eta, \sqrt{u\eta\eta'}\,\rangle$ for some unit $u \in k_1$ and $q_1 = 2$. Here $\zeta$ is a root of unity generating the torsion subgroup $W_L$ of $E_L$.

Next consider the case where $k_1$ is complex, and let $\varepsilon$ denote the fundamental unit of $k_2$. Then $\pm \varepsilon$ stays fundamental in $L$ by the argument above.

Let $\eta$ be a fundamental unit in $K_1$. If $\pm \eta$ became a square in $L$, then clearly $L/K_1$ could not be essentially ramified. Thus if we have $q_1 \geq 4$, then $\pm \varepsilon \eta = \alpha^2$ is a square in $L$. Applying $\tau$ to this relation we find that $-1 = \varepsilon \varepsilon'$ is a square in $L$, contradicting the assumption that $L$ does not contain $\sqrt{-1}$. □

**Proposition 7.** *Suppose that $q_2 = 1$. Then $K_2/k_2$ is essentially ramified if and only if $\kappa_2 = 1$; if $K_2/k_2$ is not essentially ramified, then $\kappa_2 = \langle [\mathfrak{b}] \rangle$, where $K_2 = k_2(\sqrt{\beta}\,)$ and $(\beta) = \mathfrak{b}^2$.*

*Proof.* First notice that if $K_2/k_2$ is not essentially ramified, then $\kappa_2 \neq 1$: in fact, in this case we have $(\beta) = \mathfrak{b}^2$, and if we had $\kappa_2 = 1$, then $\mathfrak{b}$ would have to be principal, say $\mathfrak{b} = (\gamma)$. This implies that $\beta = \varepsilon \gamma^2$ for some unit $\varepsilon \in k_2$, which in view of $q_2 = 1$ implies that $\varepsilon$ must be a square. But then $\beta$ would be a square, and this is impossible.

Conversely, suppose $\kappa_2 \neq 1$. Let $\mathfrak{a}$ be a nonprincipal ideal in $k_2$ of absolute norm $a$, and assume that $\mathfrak{a} = (\alpha)$ in $K_2$. Then $\alpha^{1-\sigma^2} = \eta$ for some unit $\eta \in E_2$, and similarly $\alpha^{\sigma - \sigma^3} = \eta'$, where $\eta'$ is a unit in $E_2'$. But then $\eta \eta' = \alpha^{1+\sigma-\sigma^2-\sigma^3} \stackrel{2}{=} N_{L/k}\alpha = \pm N_{L/k}\mathfrak{a} = \pm a^2 \stackrel{2}{=} \pm 1$ in $L^\times$, where $\stackrel{2}{=}$ means equal up to a square in $L^\times$. Thus $\pm \eta \eta'$ is a square in $L$, so our assumption that $q_2 = 1$ implies that $\pm \eta \eta'$ must be a square in $k_2$. The same argument show that $\pm \eta/\eta'$ is a square in $k_2$, hence we find $\eta \in k_2$. Thus $\alpha^{1-\sigma^2}$ is fixed by $\sigma^2$ and so $\beta := \alpha^2 \in k_2$. This gives $K_2 = k_2(\sqrt{\beta})$, hence $K_2/k_2$ is not essentially ramified, and moreover, $\mathfrak{a} \sim \mathfrak{b}$. □

From now on assume that $k$ is one of the imaginary quadratic fields of type A) or B) as explained in the Introduction. Let

$k_1 = \mathbb{Q}(\sqrt{d_1})$ and $k_2 = \mathbb{Q}(\sqrt{d_2 d_3})$ in case A), and

$k_1 = \mathbb{Q}(\sqrt{d_3})$ and $k_2 = \mathbb{Q}(\sqrt{d_1 d_2})$ in case B).

Then there exist two unramified cyclic quartic extensions of $k$ which are $D_4$ over $\mathbb{Q}$ (see Proposition 2). Let us say a few words about their construction. Consider e.g. case B); by Rédei's theory (see [12]), the $C_4$-factorization $d = d_1 d_2 \cdot d_3$ implies that unramified cyclic quartic extensions of $k = \mathbb{Q}(\sqrt{d})$ are constructed by choosing a "primitive" solution $(x, y, z)$ of $d_1 d_2 X^2 + d_3 Y^2 = Z^2$ and putting $L = k(\sqrt{d_1 d_2}, \sqrt{\alpha})$ with $\alpha = z + x\sqrt{d_1 d_2}$ (primitive here means that $\alpha$ should not be divisible by rational integers); the other unramified cyclic quartic extension is then $\widetilde{L} = k(\sqrt{d_1 d_2}, \sqrt{d_1 \alpha})$. If we put $\beta = \frac{1}{2}(z + y\sqrt{d_3})$, then it is an elementary exercise to show that $\alpha\beta$ is a square in $L$, hence we also have $L = k(\sqrt{d_3}, \sqrt{\beta})$ etc. If $d_3 = -4$, then it is easy to see that we may choose $\beta$ as the fundamental unit of $k_2$; if $d_3 \ne -4$, then genus theory says that a) the class number $h$ of $k_2$ is twice an odd number $u$; and b) the prime ideal $\mathfrak{p}_3$ above $d_3$ in $k_2$ is in the principal genus, so $\mathfrak{p}_3^u = (\pi_3)$ is principal. Again it can be checked that $\beta = \pm \pi_3$ for a suitable choice of the sign.

**Example.** Consider the case $d = -31 \cdot 5 \cdot 8$; here $\pi_3 = \pm(3 + 2\sqrt{10})$, and the positive sign is correct since $3 + 2\sqrt{10} \equiv (1 + \sqrt{10})^2 \mod 4$ is primary. The minimal polynomial of $\sqrt{\pi_3}$ is $f(x) = x^4 - 6x^2 - 31$: compare Table 1.

The fields $K_2 = k_2(\sqrt{\alpha})$ and $\widetilde{K}_2 = k_2(\sqrt{d_2 \alpha})$ will play a dominant role in the proof below; they are both contained in $M = F(\sqrt{\alpha})$ for $F = k_2(\sqrt{d_2})$, and it is the ambiguous class group $\mathrm{Am}(M/F)$ that contains the information we are interested in.

**Lemma 6.** *The field $F$ has odd class number (even in the strict sense), and we have $\# \mathrm{Am}(M/F) \mid 2$. In particular, $\mathrm{Cl}_2(M)$ is cyclic (though possibly trivial).*

*Proof.* The class group in the strict sense of $k_2$ is cyclic of order 2 by Rédei's theory [12] (since $(d_2/p_3) = (d_3/p_2) = -1$ in case A) and $(d_1/p_2) = (d_2/p_1) = -1$ in case B)). Since $F$ is the Hilbert class field of $k_2$ in the strict sense, its class number in the strict sense is odd.

Next we apply the ambiguous class number formula. In case A), $F$ is complex, and exactly the two primes above $d_3$ ramify in $M/F$. Note that $M = F(\sqrt{\alpha})$ with $\alpha$ primary of norm $d_3 y^2$; there are four primes above $d_3$ in $F$, and exactly two of them divide $\alpha$ to an odd power, so $t = 2$ by the decomposition law in quadratic Kummer extensions. By Proposition 4 and the remarks following it, $\# \mathrm{Am}_2(M/F) = 2/(E:H) \le 2$, and $\mathrm{Cl}_2(M)$ is cyclic.

In case B), however, $F$ is real; since $\alpha \in k_2$ has norm $d_3 y^2 < 0$, it has mixed signature, hence there are exactly two infinite primes that ramify in $M/F$. As in case A), there are two finite primes above $d_3$ that ramify in $M/F$, so we get $\# \mathrm{Am}_2(M/F) = 8/(E:H)$. Since $F$ has odd class number in the strict sense, $F$ has units of independent signs. This implies that the group of units that are positive at the two ramified infinite primes has $\mathbb{Z}$-rank 2, i.e. $(E:H) \geq 4$ by consideration of the infinite primes alone. In particular, $\#\operatorname{Am}_2(M/F) \leq 2$ in case B). □

Next we derive some relations between the class groups of $K_2$ and $\widetilde{K}_2$; these relations will allow us to use each of them as our field $K$ in Theorem 1.

**Proposition 8.** *Let $L$ and $\widetilde{L}$ be the two unramified cyclic quartic extensions of $k$, and let $K_2$ and $\widetilde{K}_2$ be two quadratic extensions of $k_2$ in $L$ and $\widetilde{L}$, respectively, which are not normal over $\mathbb{Q}$.*

  a) *We have $4 \mid h(K_2)$ if and only if $4 \mid h(\widetilde{K}_2)$;*
  b) *If $4 \mid h(K_2)$, then one of $\operatorname{Cl}_2(K_2)$ or $\operatorname{Cl}_2(\widetilde{K}_2)$ has type $(2,2)$, whereas the other is cyclic of order $\geq 4$.*

*Proof.* Notice that the prime dividing $\operatorname{disc}(k_1)$ splits in $k_2$. Throughout this proof, let $\mathfrak{p}$ be one of the primes of $k_2$ dividing $\operatorname{disc}(k_1)$.

If we write $K_2 = k_2(\sqrt{\alpha})$ for some $\alpha \in k_2$, then $\widetilde{K}_2 = k_2(\sqrt{d_2\alpha})$. In fact, $K_2$ and $\widetilde{K}_2$ are the only extensions $F/k_2$ of $k_2$ with the properties
1. $F/k_2$ is a quadratic extension unramified outside $\mathfrak{p}$;
2. $kF/k$ is a cyclic extension.
Therefore it suffices to observe that if $k_2(\sqrt{\alpha})$ has these properties, then so does $k_2(\sqrt{d_2\alpha})$. But this is elementary.

In particular, the compositum $M = K_2\widetilde{K}_2 = k_2(\sqrt{d_2}, \sqrt{\alpha})$ is an extension of type $(2,2)$ over $k_2$ with subextensions $K_2$, $\widetilde{K}_2$ and $F = k_2(\sqrt{d_2})$. Clearly $F$ is the unramified quadratic extension of $k_2$, so both $M/K_2$ and $M/\widetilde{K}_2$ are unramified. If $K_2$ had 2-class number 2, then $M$ would have odd class number, and $M$ would also be the 2-class field of $\widetilde{K}_2$. Thus $2 \parallel h(K_2)$ implies that $2 \parallel h(\widetilde{K}_2)$. This proves part a) of the proposition.

Before we go on, we give a Hasse diagram for the fields occurring in this proof:

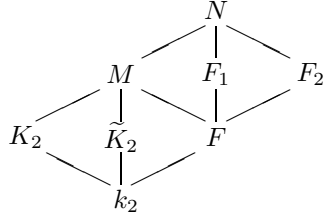

Now assume that $4 \mid h(K_2)$. Since $\operatorname{Cl}_2(M)$ is cyclic by Lemma 6, there is a unique quadratic unramified extension $N/M$, and the uniqueness implies at once that $N/k_2$ is normal. Hence $G = \operatorname{Gal}(N/k_2)$ is a group of order 8 containing a subgroup of type $(2,2) \simeq \operatorname{Gal}(N/F)$: in fact, if $\operatorname{Gal}(N/F)$ were cyclic, then the primes ramifying in $M/F$ would also ramify in $N/M$ contradicting the fact that $N/M$ is unramified. There are three groups satisfying these conditions: $G = (2,4)$, $G = (2,2,2)$ and $G = D_4$. We claim that $G$ is non-abelian; once we have proved this, it follows that exactly one of the groups $\operatorname{Gal}(N/K_2)$ and $\operatorname{Gal}(N/\widetilde{K}_2)$ is cyclic, and that the other is not, which is what we want to prove.

So assume that $G$ is abelian. Then $M/F$ is ramified at two finite primes $\mathfrak{q}$ and $\mathfrak{q}'$ of $F$ dividing $\mathfrak{p}$ (in $k_2$); if $F_1$ and $F_2$ denote the quadratic subextensions of $N/F$ different from $M$ then $F_1/F$ and $F_2/F$ must be ramified at a finite prime (since $F$ has odd class number in the strict sense: see Lemma 6); since both $F_1$ and $F_2$ are normal (even abelian) over $k_2$, ramification at $\mathfrak{q}$ implies ramification at the conjugated ideal $\mathfrak{q}'$. Hence both $\mathfrak{q}$ and $\mathfrak{q}'$ ramify in $F_1/F$ and $F_2/F$, and since they also ramify in $M/F$, they must ramify completely in $N/F$, again contradicting the fact that $N/M$ is unramified.

We have proved that $\mathrm{Cl}_2(K_2)$ and $\mathrm{Cl}_2(\widetilde{K}_2)$ contain subgroups of type (4) and $(2,2)$, respectively. Now we wish to apply Proposition 5. But we have to compute $\#\,\mathrm{Am}_2(\widetilde{K}_2/k_2)$. Since the class number of $\widetilde{K}_2$ is even, it is sufficient to show that $\#\,\mathrm{Am}_2(\widetilde{K}_2/k_2) \leq 2$. In case A), there is exactly one ramified prime (it divides $d_1$), hence $\#\,\mathrm{Am}_2(\widetilde{K}_2/k_2) = 2/(E:H) \leq 2$. In case B), there are two ramified primes (one is infinite, the other divides $d_3$), hence $\#\,\mathrm{Am}_2(\widetilde{K}_2/k_2) = 4/(E:H)$; but $-1$ is not a norm residue at the ramified infinite prime, hence $(E:H) \geq 2$ and $\#\,\mathrm{Am}_2(\widetilde{K}_2/k_2) \leq 2$ as claimed.

Now Proposition 5 implies that $\mathrm{Cl}_2(K_2)$ is cyclic of order $\geq 4$, and that $\mathrm{Cl}_2(\widetilde{K}_2) \simeq (2,2)$. This concludes our proof. □

**Proposition 9.** *Assume that $k$ is one of the imaginary quadratic fields of type A) or B) as explained in the Introduction. Then there exist two unramified cyclic quartic extensions of $k$. Let $L$ be one of them, and write*

$k_1 = \mathbb{Q}(\sqrt{d_1})$ *and* $k_2 = \mathbb{Q}(\sqrt{d_2 d_3})$ *in case A), and*

$k_1 = \mathbb{Q}(\sqrt{d_3})$ *and* $k_2 = \mathbb{Q}(\sqrt{d_1 d_2})$ *in case B).*

*Then $h_2(L) = \frac{1}{4} h_2(k) h_2(K_1) h_2(K_2)$ unless possibly when $d_3 = -4$ in case B).*

*Proof.* Observe that $v = 0$ in case A) and B); Kuroda's class number formulas for $L/k_1$ and $L/k_2$ gives

$$h_2(L) = \frac{q_1 h_2(K_1)^2 h_2(K)}{2 h_2(k_1)^2} = \frac{q_2 h_2(K_2)^2 h_2(K)}{4 h_2(k_2)^2}$$

in case A) and

$$h_2(L) = \frac{q_1 h_2(K_1)^2 h_2(K)}{4 h_2(k_1)^2} = \frac{q_2 h_2(K_2)^2 h_2(K)}{2 h_2(k_2)^2}$$

in case B). Multiplying them together and plugging in the class number formula for $K/\mathbb{Q}$ yields

$$h_2(L)^2 = \frac{q_1 q_2}{8} \frac{h_2(K_1)^2 h_2(K_2)^2 h_2(k)^2}{h_2(k_1)^2 h_2(k_2)^2}.$$

Now $h_2(k_1) = 1$, $h_2(k_2) = 2$ and $q_1 q_2 = 2$ (by Proposition 6), and taking the square root we find $h_2(L) = \frac{1}{4} h_2(k) h_2(K_1) h_2(K_2)$ as claimed. □

## 5. Classification

In this section we apply the results obtained in the last few sections to give a proof for Theorem 1.

*Proof of Theorem 1.* Let $L$ be one of the two cyclic quartic unramified extensions of $k$, and let $N$ be the subgroup of $\mathrm{Gal}(k^2/k)$ fixing $L$. Then $N$ satisfies the assumptions of Proposition 1, thus there are only the following possibilities:

| $d(G')$ | $h_2(L)$ | $h_2(K_1)h_2(K_2)$ |
|---------|----------|---------------------|
| 1 | $2^m$ | 2 |
| 2 | $2^{m+1}$ | 4 |
| $\geq 3$ | $\geq 2^{m+2}$ | $\geq 8$ |

Here, the first two columns follow from Proposition 1, the last (which we do not claim to hold if $d_3 = -4$ in case B)) is a consequence of the class number formula of Proposition 9. In particular, we have $d(G') \geq 3$ if one of the class numbers $h_2(K_1)$ or $h_2(K_2)$ is at least 8. Therefore it suffices to examine the cases $h_2(K_2) = 2$ and $h_2(K_2) = 4$ (recall from above that $h_2(K_2)$ is always even).

We start by considering case A); it is sufficient to show that $h_2(K_1)h_2(K_2) \neq 4$. We now apply Proposition 5; notice that we may do so by the proof of Proposition 8.
a) If $h_2(K_2) = 2$, then $\#\kappa_2 = 2$ by Proposition 5, hence $q_2 = 2$ by Proposition 7 and then $q_1 = 1$ by Proposition 6. The class number formulas in the proof of Proposition 9 now give $h_2(K_1) = 1$ and $h_2(L) = 2^m$.

It can be shown using the ambiguous class number formula that $\text{Cl}_2(K_1)$ is trivial if and only if $\varepsilon_1$ is a quadratic nonresidue modulo the prime ideal over $d_2$ in $k_1$; by Scholz's reciprocity law, this is equivalent to $(d_1/d_2)_4(d_2/d_1)_4 = 1$, and this agrees with the criterion given in [1].
b) If $h_2(K_2) = 4$, we may assume that $\text{Cl}_2(K_2) = (4)$ from Proposition 8.b). Then $\#\kappa_2 = 2$ by Proposition 5, $q_2 = 2$ by Proposition 7 and $q_1 = 1$ by Proposition 6. Using the class number formula we get $h_2(K_1) = 2$ and $h_2(L) = 2^{m+2}$.

Thus in both cases we have $h_2(K_1)h_2(K_2) \neq 4$, and by the table at the beginning of this proof this implies that $\text{rank}\,\text{Cl}_2(k^1) \neq 2$ in case A).

Next we consider case B); here we have to distinguish between $d_3 \neq -4$ (case $B_1$) and $d_3 = -4$ (case $B_2$).

Let us start with case $B_1$).
a) If $h_2(K_2) = 2$, then $\#\kappa_2 = 2$, $q_2 = 2$ and $q_1 = 1$ as above. The class number formula gives $h_2(K_1) = 2$ and $h_2(L) = 2^{m+1}$.
b) If $\text{Cl}_2(K_2) = (4)$ (which we may assume without loss of generality by Proposition 8.b)) then $\#\kappa_2 = 2$, $q_2 = 2$ and $q_1 = 1$, again exactly as above. This implies $h_2(K_1) = 4$ and $h_2(L) = 2^{m+3}$.

Finally, consider case $B_2$).
Here we apply Kuroda's class number formula (see [10]) to $L/k_1$, and since $h_2(k_1) = 1$ and $h_2(K_1) = h_2(K_1')$, we get $h_2(L) = \frac{1}{2}q_1 h_2(K_1)^2 h_2(k) = 2^m q_1 h_2(K_1)^2$. From $K_2 = k_2(\sqrt{\varepsilon})$ (for a suitable choice of $L$; the other possibility is $\widetilde{K}_2 = k_2(\sqrt{d_2 \varepsilon})$), where $\varepsilon$ is the fundamental unit of $k_2$, we deduce that the unit $\varepsilon$, which still is fundamental in $k$, becomes a square in $L$, and this implies that $q_1 \geq 2$. Moreover, we have $K_1 = k_1(\sqrt{\pi \lambda})$, where $\pi, \lambda \equiv 1 \bmod 4$ are prime factors of $d_1$ and $d_2$ in $k_1 = \mathbb{Q}(i)$, respectively. This shows that $K_1$ has even class number, because $K_1(\sqrt{\pi})/K_1$ is easily seen to be unramified.

Thus $2 \mid q_1$, $2 \mid h_2(K_1)$, and so we find that $h_2(L)$ is divisible by $2^m \cdot 2 \cdot 4 = 2^{m+3}$. In particular, we always have $d(G') \geq 3$ in this case.

This concludes the proof. □

The referee (whom we'd like to thank for a couple of helpful remarks) asked whether $h_2(K) = 2$ and $h_2(K) > 2$ infinitely often. Let us show how to prove that both possibilities occur with equal density.

Before we can do this, we have to study the quadratic extensions $K_1$ and $\widetilde{K}_1$ of $k_1$ more closely. We assume that $d_2 = p$ and $d_3 = r$ are odd primes in the following, and then say how to modify the arguments in the case $d_2 = 8$ or $d_3 = -8$. The primes $p$ and $r$ split in $k_1$ as $p\mathcal{O}_1 = \mathfrak{p}\mathfrak{p}'$ and $r\mathcal{O}_1 = \mathfrak{r}\mathfrak{r}'$. Let $h$ denote the odd class number of $k_1$ and write $\mathfrak{p}^h = (\pi)$ and $\mathfrak{r}^h = (\rho)$ for primary elements $\pi$ and $\rho$ (this is can easily be proved directly, but it is also a very special case of Hilbert's first supplementary law for quadratic reciprocity in fields $K$ with odd class number $h$ (see [7]): if $\mathfrak{a}^h = \alpha\mathcal{O}_K$ for an ideal $\mathfrak{p}$ with odd norm, then $\alpha$ can be chosen primary (i.e. congruent to a square mod $4\mathcal{O}_K$) if and only if $\mathfrak{a}$ is primary (i.e. $[\varepsilon/\mathfrak{a}] = +1$ for all units $\varepsilon \in \mathcal{O}_K^\times$, where $[\cdot/\cdot]$ denotes the quadratic residue symbol in $K$)). Let $[\cdot/\cdot]$ denote the quadratic residue symbol in $k_1$. Then $[\pi/\rho][\pi'/\rho] = [p/\rho] = (p/r) = -1$, so we may choose the conjugates in such a way that $[\pi/\rho] = +1$ and $[\pi'/\rho] = [\pi/\rho'] = -1$.

Put $K_1 = k_1(\sqrt{\pi\rho})$ and $\widetilde{K}_1 = k_1(\sqrt{\pi\rho'})$; we claim that $h_2(\widetilde{K}_1) = 2$. This is equivalent to $h_2(\widetilde{L}_1) = 1$, where $\widetilde{L}_1 = k_1(\sqrt{\pi}, \sqrt{\rho'})$ is a quadratic unramified extension of $\widetilde{K}_1$. Put $\widetilde{F}_1 = k_1(\sqrt{\pi})$ and apply the ambiguous class number formula to $\widetilde{F}_1/k_1$ and $\widetilde{L}_1/\widetilde{F}_1$: since there is only one ramified prime in each of these two extensions, we find $\mathrm{Am}(\widetilde{F}_1/k_1) = \mathrm{Am}(\widetilde{L}_1/\widetilde{F}_1) = 1$; note that we have used the assumption that $[\pi/\rho'] = -1$ in deducing that $\mathfrak{r}'$ is inert in $\widetilde{F}_1/k_1$.

In our proof of Theorem 1 we have seen that there are the following possibilities when $h_2(K_2) \mid 4$:

| $q_2$ | $\mathrm{Cl}_2(K_2)$ | $q_1$ | $h_2(K_1)$ | $\widetilde{q}_2$ | $\mathrm{Cl}_2(\widetilde{K}_2)$ | $h_2(L)$ |
|---|---|---|---|---|---|---|
| 2 | (2) | 1 | 2 | 2 | (2) | $2^{m+1}$ |
| 2 | (4) | 1 | 4 | ? | (2,2) | $2^{m+3}$ |

In order to decide whether $\widetilde{q}_2 = 1$ or $\widetilde{q}_2 = 2$, recall that we have $h_2(K_1) = 4$; thus $\widetilde{K}_1$ must be the field with 2-class number 2, and this implies $h_2(\widetilde{L}) = 2^{m+2}$ and $\widetilde{q}_2 = 1$. In particular we see that $4 \mid h_2(K_2)$ if and only if $4 \mid h_2(K_1)$ as long as $K_1 = k_1(\sqrt{\pi\rho})$ with $[\pi/\rho] = +1$.

The ambiguous class number formula shows that $\mathrm{Cl}_2(K_1)$ is cyclic, thus $4 \mid h_2(K_1)$ if and only if $2 \mid h_2(L_1)$, where $L_1 = K_1(\sqrt{\pi})$ is the quadratic unramified extension of $K_1$. Applying the ambiguous class number formula to $L_1/F_1$, where $F_1 = k_1(\sqrt{\pi})$, we see that $2 \mid h_2(L_1)$ if and only if $(E : H) = 1$. Now $E$ is generated by a root of unity (which always is a norm residue at primes dividing $r \equiv 1 \bmod 4$) and a fundamental unit $\varepsilon$. Therefore $(E : H) = 1$ if and only if $\{\varepsilon/\mathfrak{R}_1\} = \{\varepsilon/\mathfrak{R}_2\} = +1$, where $\mathfrak{r}\mathcal{O}_{F_1} = \mathfrak{R}_1\mathfrak{R}_2$ and where $\{\cdot/\cdot\}$ denotes the quadratic residue symbol in $F_1$. Since $\{\varepsilon/\mathfrak{R}_1\}\{\varepsilon/\mathfrak{R}_2\} = [\varepsilon/\mathfrak{r}] = +1$, we have proved that $4 \mid h_2(K_1)$ if and only if the prime ideal $\mathfrak{R}_1$ above $\mathfrak{r}$ splits in the quadratic extension $F_1(\sqrt{\varepsilon})$. But if we fix $p$ and $q$, this happens for exactly half of the values of $r$ satisfying $(p/r) = -1$, $(q/r) = +1$.

If $d_2 = 8$ and $p = 2$, then $2\mathcal{O}_{k_1} = \mathfrak{z}\mathfrak{z}'$, and we have to choose $\mathfrak{z}^h = (\pi)$ in such a way that $k_1(\sqrt{\pi})/k_1$ is unramified outside $\mathfrak{p}$. The residue symbols $[\alpha/\mathfrak{z}]$ are defined as Kronecker symbols via the splitting of $\mathfrak{z}$ in the quadratic extension $k_1(\sqrt{\alpha})/k_1$. With these modifications, the above arguments remain valid.

(E. Benjamin) Mathematics Department, Unity College
*E-mail address*: `benjamin@mint.net`

(F. Lemmermeyer) Erwin-Rohde-Str. 19, Heidelberg, Germany
*E-mail address*: `hb3@ix.urz.uni-heidelberg.de`

(C. Snyder) Department of Mathematics and Statistics, University of Maine, and, Research Institute of Mathematics, Orono
*E-mail address*: `snyder@gauss.umemat.maine.edu`